\documentclass{amsart}
\usepackage[nosections]{macs-fmw}
\allowdisplaybreaks[2]

\date{20 May 2002}

\begin{document}

\title{Continuous-Trace Groupoid Crossed Products}

\author{Igor Fulman}
\address{Department of Mathematics \\Arizona State
  University\\Tempe\\Arizona 85287-1804}
\email{ifulman@math.la.asu.edu}
\author{Paul S. Muhly}
\address{Department of Mathematics\\University of Iowa\\Iowa City, IA
  52242}
\email{muhly@math.uiowa.edu}
\author{Dana P. Williams}
\address{Department of Mathematics\\Dartmouth College\\Hanover, NH
  03755-3551} 
\email{dana.williams@dartmouth.edu}

\begin{abstract}
  Let $G$ be a second countable, locally compact groupoid with Haar
  system, and let $\mathcal{A}$ be a bundle of $C^{\ast}$-algebras
  defined over the unit space of $G$ on which $G$ acts continuously.
  We determine conditions under which the associated crossed product
  $C^{\ast}(G;\mathcal{A})$ is a continuous trace $C^{\ast}$-algebra.
\end{abstract}

\maketitle

\section{Introduction}

Throughout this note, $G$ will denote a second countable, locally
compact groupoid with Haar system $\{\lambda^{u}\}_{u\in G^{(0)}}$.
Also, we shall fix a bundle\footnote{We follow
  \cite{fd:representations} for the general theory of Banach- and
  $C^{\ast}$-bundles. However, we adopt the increasingly popular
  convention that bundles are to be denoted by calligraphic letters.
  The fibres in a given bundle are then denoted by the corresponding
  roman letter. Thus if $\mathcal{A}$ is a bundle of
  $C^{\ast}$-algebras, say, over a space $X$, then the fibre over
  $x\in X$ will be denoted $A(x)$.} $\mathcal{A}$ of (separable)
$C^{\ast}$-algebras over the unit space $G^{(0)}$ of $G$. We shall
write $p$ for the projection of $\mathcal{A}$ onto $G^{(0)}$. We shall
assume that there is a continuous action, denoted $\sigma$, of $G$ on
$\mathcal{A}$. This means the following: First of all, $\sigma$ must
be a homomorphism from $G$ into the isomorphism groupoid of
$\mathcal{A}$, $\Iso(\mathcal{A})$, so that, in particular,
$\sigma_{\gamma}:A(s(\gamma))\rightarrow A(r(\gamma))$ is a
$C^{\ast}$-isomorphism for each $\gamma\in G$. Secondly, let $s^{\ast
}(\mathcal{A})$ and $r^{\ast}(\mathcal{A})$ be the bundles on $G$
obtained by
pulling back $\mathcal{A}$ via $s$ and $r$, so that $s^{\ast}(\mathcal{A}%
)=\{(\gamma,a)\mid a\in A(s(\gamma))\}$ and similarly for $r^{\ast
}(\mathcal{A})$. Then $\sigma$ determines a bundle map
$\sigma^{\ast}:s^{\ast }(\mathcal{A})\rightarrow
r^{\ast}(\mathcal{A})$ by the formula $\sigma^{\ast
}(\gamma,a)=(\gamma,\sigma_{\gamma}(a))$. The continuity assumption
that we make is that for each continuous section $f$ of
$s^{\ast}(\mathcal{A})$, $\sigma^{\ast}\circ f$ is a continuous
section of $r^{\ast}(\mathcal{A})$.

Let $C_{c}(G,r^{\ast}(\mathcal{A}))$ denote the space of continuous
sections
of $r^{\ast}(\mathcal{A})$ with compact support and for $f,g\in C_{c}%
(G,r^{\ast}(\mathcal{A}))$, set%
\[
f\ast g(\gamma):=\int
f(\eta)\sigma_{\eta}(g(\eta^{-1}\gamma))\,d\lambda ^{r(\gamma)}(\eta)
\]
and set%
\[
f^{\ast}(\gamma)=\sigma_{\gamma}(f(\gamma^{-1})^{\ast})\text{.}%
\]
Then,  with respect  to these  operations and  pointwise  addition and
scalar  multiplication,   $C_{c}(G,r^{\ast}(\mathcal{A}))$  becomes  a
topological $\ast  $-algebra in the inductive limit  topology to which
Renault's   disintegration  theorem   \cite{ren:jot87}   applies.  The
enveloping  $C^{\ast}$-algebra of  $C_{c}(G,r^{\ast}(\mathcal{A}))$ is
called   the   \emph{crossed   product   of}  $G$   \emph{acting   on}
$\mathcal{A}$     and     is     denoted     $C^{\ast}(G;\mathcal{A})$
\cite{ren:jot87}. The basic problem we study in this note is

\begin{question}
\label{FundQuest}Under what circumstances on $G$ and $\mathcal{A}$ is
$C^{\ast}(G;\mathcal{A})$ a continuous trace $C^{\ast}$-algebra?
\end{question}

The first systematic investigation into questions of this type that we
know of is Green's pioneering study \cite{gre:pjm77}. There,
Green deals with the case where $G$ is the transformation group
groupoid obtained by letting a locally compact \emph{group} $H$ act on
a locally compact Hausdorff space $X$ and where $\mathcal{A}$ is the
trivial line bundle over $G^{(0)}=X$. His principal result asserts
that if $H$ acts freely on $X$, then $C^{\ast}(G;\mathcal{A})$ is
continuous trace if and only if the action of $H$ is proper. In this
event, in fact, $C^{\ast}(G;\mathcal{A})$ is strongly Morita
equivalent to $C_{0}(X/H)$.

Another precedent to this investigation is the work of the second two
authors in \cite{muhwil:ms90}. Here, the hypothesis is that $G$ is a
\emph{principal} groupoid and the bundle $\mathcal{A}$ is again the
trivial line bundle. It was shown that
$C^{\ast}(G;\mathcal{A})=C^{\ast}(G)$ has continuous trace if and only
if the action of $G$ on $G^{(0)}$ is proper. In this event, again,
$C^{\ast}(G)$ is strongly Morita equivalent to $C_{0}(G^{(0)}/G)$.

The first example where the bundle $\mathcal{A}$ is nontrivial was
considered by Raeburn and Rosenberg in \cite{raeros:tams88}. They
considered a locally compact \emph{group} acting on a continuous trace
$C^{\ast}$-algebra $A$ and showed that if the natural action of $G$ on
the spectrum of $A$, $\hat{A}$, is free and proper, then the cross
product $C^{\ast}$-algebra $A\rtimes G$ has continuous trace. In
\cite{olerae:jfa90}, Olesen and Raeburn proved a conditioned converse:
if the group $G$ is abelian and acts freely on $\hat{A}$, then if
$A\rtimes G$ is continuous trace, the action of $G$ on $\hat{A}$ must
be proper.  Quite recently, Deicke \cite{dei:jot00} used non-abelian
duality theory to remove the hypothesis that $G$ is abelian. Thus, the
best result in this direction is: If $G$ acts on $A$ yielding a free
action on $\hat{A}$, then $A\rtimes G$ is continuous trace if and only
if the action of $G$ on $\hat{A}$ is proper.

Our objective in this note, Theorem~\ref{Main}, is to prove a result
that contains all of these examples as special cases --- and much
more, as well. It is based on two hypotheses and some ancillary
considerations that we will elaborate. The first hypothesis is

\begin{hypothesis}
\label{Hyp1} The $C^{\ast}$-algebra $C_{0}(G^{(0)},\mathcal{A})$ is
continuous trace.
\end{hypothesis}

Here, $C_{0}(G^{(0)},\mathcal{A})$ denotes the $C^{\ast}$-algebra of
continuous sections of the bundle $\mathcal{A}$ that vanish at
infinity on $G^{(0)}$. An hypothesis on the bundle $\mathcal{A}$ of
this nature is natural and reasonable, in view of the fact that in the
trivial case $G=G^{(0)}$, we have $C^{*}(G;\mathcal{A}) \cong
C_{0}(G^{(0)},\mathcal{A})$.  Furthermore, we note that a
compact group can act on an antiliminal $C^{\ast}$-algebra in such a
way that the crossed product is continuous trace \cite{tak:jfa71}.

To state our second hypothesis, we need a couple of remarks about the
spectrum of $C_{0}(G^{(0)},\mathcal{A})$. We shall denote it by $X$
throughout this note. Observe that $X$ may be expressed as the
disjoint union of spectra $\coprod_{u\in G^{(0)}}\specnp{A(u)}$ and the natural
projection $\hat{p}$ from $X$ to $G^{(0)}$ is continuous and open
\cite{lee:iumj76,nil:iumj96}.  The groupoid $G$ acts on $X$ (using the
map $\hat{p}$) as follows. If $x\in X$, we shall write $x=[\pi_{x}]$
in order to specify a particular irreducible representation in the
equivalence class represented by $x$. Then if $X\ast G$ denotes the
space $\{(x,\gamma)\in X\times G\mid\hat {p}(x)=r(x)\}$ and if
$(x,\gamma)\in X\ast G$, then $x\cdot\gamma$ is defined to be
$[\pi_{x}\circ\sigma_{\gamma}]$. The fact that the action of $G$ on
$X$ is well defined and is continuous is easily checked (as in
\cite[Lemma~7.1]{rw:morita} for example).  Our second hypothesis~is

\begin{hypothesis}
\label{Hyp2}The action of $G$ on $X$ is free.
\end{hypothesis}

For $x\in X$, we shall write $\mathcal{K}(x)$ for the quotient $\mathcal{A}%
(\hat{p}(x))/\ker(\pi_{x})$. Then $\mathcal{K}(x)$ is well defined
(i.e., it is
independent of the choice of $\pi_{x}$) and is an elementary $C^{\ast}%
$-algebra. Hypothesis~\ref{Hyp1} guarantees that the collection $\{\mathcal{K}%
(x)\}_{x\in X}$ may be given the structure of an elementary
$C^{\ast}$-algebra bundle over $X$ satisfying Fell's condition
\cite[Proposition~10.5.8]{dix:cs-algebras}.  As a result, we find that
$C_{0}(G^{(0)},\mathcal{A})$ is naturally isomorphic to
$C_{0}(X,\mathcal{K})$. The action of $G$ on $X$ induces one on
$\mathcal{K}$ that we shall use. We find it preferable to express this
in terms of the \emph{action groupoid} $X\ast G$ defined as follows.

As a set, $X\ast G:=\{(x,\gamma)\mid\hat{p}(x)=r(\gamma)\}$ and the
groupoid operations are defined by the formulae
\begin{align*}
  (x,\alpha)(x\alpha,\beta)  &  :=(x,\alpha\beta)\;\text{and}\\
(x,\alpha)^{-1}  &  :=(x\alpha,\alpha^{-1})\text{,}%
\end{align*}
$(x,\alpha)$, $(x,\beta)\in X\ast G$. Note, in particular, that the
unit space
of $X\ast G$ may be identified with $X$ via: $(x,\hat{p}%
(x))\longleftrightarrow x$. Note, too, that the range and source maps
on $X\ast G$, denoted $\tilde{r}$ and $\tilde{s}$, are given by the
equations $\tilde{r}(x,\alpha)=(x,r(\alpha))$ and
$\tilde{s}(x,\alpha)=(x\alpha ,s(\alpha))$. The groupoid $X\ast G$ in
the product topology is clearly locally compact, Hausdorff, and
separable. It has a Haar system $\{\tilde {\lambda}^{x}\}_{x\in X}$
given by the formula $\tilde{\lambda}^{x}=\delta
_{x}\times\lambda^{\hat{p}(x)}$. Observe that the action of $G$ on $X$
is free (resp. proper) iff $X\ast G$ is principal (resp. proper).

The groupoid $X\ast G$ acts on $\mathcal{K}$ via the formula%
\[
\tilde{\sigma}_{(x,\gamma)}(k):=\sigma_{\gamma}(a)+\ker(\pi_{x})\text{,}%
\]
where $k=a+\ker\pi_{x\gamma}$ lies in $K(\tilde{s}(x,\gamma))$. Note
that this action is well defined since
$\pi_{x}\circ\sigma_{\gamma}=\pi_{x\cdot\gamma}$.  We promote this
action of $X\ast G$ on $\mathcal{K}$ to one on $X\ast
\mathcal{K}:=\{(x,k)\mid k\in K(x)\}$: $(x,k)\cdot(x,\gamma):=(x\cdot
\gamma,\tilde{\sigma}_{(x,\gamma)}^{-1}(k))$. If the action of $G$ on
$X$ is free and proper, then the action of $X\ast G$ on
$X\ast\mathcal{K}$ is also free and proper. In this case, we write
$\mathcal{K}^{X}$ for the quotient space $X\ast\mathcal{K}/X\ast G$.
Then $\mathcal{K}^{X}$ is naturally a bundle of elementary
$C^{\ast}$-algebras over $X/G$. In fact, using the methods of
Theorem~1.1 of \cite{raeros:tams88}, it is easy to see that
$\mathcal{K}^{X}$ satisfies Fell's condition. Thus, in particular,
$C_{0}(X/G,\mathcal{K}^{X})$ has continuous trace, if $G$ acts on $X$
freely and properly.

With these preliminaries at our disposal, we are able to state the
main result of this paper as

\begin{thm}
\label{Main}Under Hypotheses \ref{Hyp1}~and \ref{Hyp2}, $C^{\ast
}(G;\mathcal{A})$ has continuous trace if and only if the action of $G$
on $X$ is proper. In this event, $C^{\ast}(G;\mathcal{A})$ is strongly
Morita equivalent to $C_{0}(X/G,\mathcal{K}^{X})$, where
$\mathcal{K}^{X}$ is the elementary $C^{\ast}$-bundle over $X/G,$
satisfying Fell's condition, that was just defined.
\end{thm}

\section{Sufficiency in Theorem~\ref{Main}}

First we reduce the proof of Theorem~\ref{Main} to the case when $G$
is a principal groupoid. This reduction is accomplished
with the aid of

\begin{thm}
\label{Reduction}In the notation established above, $C^{\ast}(G;\mathcal{A})$
is isomorphic to $C^{\ast}(X\ast G;\mathcal{K})$.
\end{thm}

\begin{proof}
  For $u\in G^{(0)}$, identify $A(u)$ with
  $C_{0}(\hat{p}^{-1}(u),\mathcal{K})$.  Also, set
  $C_{cc}(G,r^{\ast}(\mathcal{A})):=\{f\in C_{c}(G,r^{\ast
  }(\mathcal{A}))\mid(x,\gamma)\rightarrow\left\| f(\gamma)(x)\right\|
  $ has
compact support in $X\ast G\}$. Recall that $r^{\ast}(\mathcal{A}%
)=\{(\gamma,a)\mid a\in A(r(\gamma))\}$ and, likewise,
$\tilde{r}^{\ast }(\mathcal{K})=\{((x,\gamma),a)\mid a\in K(x)\}$. So,
in the definition of $C_{cc}(G,r^{\ast}(\mathcal{A}))$, $f(\gamma)$,
which nominally is in $A(r(\gamma))$, is to be viewed in
$C_{0}(\hat{p}^{-1}(u),\mathcal{K})$. In particular, $f(\gamma)(x)$
lies in $K(x)$, when $\hat{p}(x)=r(\gamma)$. Now it is straightforward
to check that $C_{cc}(G,r^{\ast}(\mathcal{A}))$ is a
subalgebra of $C_{c}(G,r^{\ast}(\mathcal{A}))$ that is dense in $C_{c}%
(G,r^{\ast}(\mathcal{A}))$ in the inductive limit topology. Further,
if we define $\Psi:C_{c}(X\ast
G,\tilde{r}^{\ast}(\mathcal{K}))\rightarrow
C_{cc}(G,r^{\ast}(\mathcal{A}))$ by the formula%
\[
\lbrack\Psi(f)(\gamma)](x)=f(x,\gamma)\text{,}%
\]
$f\in C_{c}(X\ast G,\tilde{r}^{\ast}(\mathcal{K}))$, then $\Psi$ is an
algebra $\ast$-homomorphism that is continuous in the inductive limit
topologies on $C_{c}(X\ast G,\tilde{r}^{\ast}(\mathcal{K}))$ and
$C_{cc}(G,r^{\ast }(\mathcal{A}))$. Hence by Renault's disintegration
theorem \cite[Theorem 4.1]{ren:jot87}, $\Psi$ extends to a
$C^{\ast}$-homomorphism from $C^{\ast}(X\ast G;\mathcal{K})$ into
$C^{\ast}(G;\mathcal{A})$.

The inverse $\Phi:C_{cc}(G,r^{\ast}(\mathcal{A}))\rightarrow
C_{c}(X\ast
G,\tilde{r}^{\ast}(\mathcal{K}))$ to $\Psi$ is given formally by the formula%
\[
\Phi(f)(x,\gamma)=f(\gamma)(x)\text{,}%
\]
$f\in C_{cc}(G,r^{\ast}(\mathcal{A}))$; i.e., at the level of $C_{c}$,
$\Phi\circ\Psi$ and $\Psi\circ\Phi$ are the identity maps on the
appropriate algebras. The problem is that \emph{a priori} $\Phi$ is
not continuous in the inductive limit topologies. However, $\Phi$
\emph{is} manifestly continuous
with respect to the so-called $L^{I}$-norms on $C_{cc}(G,r^{\ast}%
(\mathcal{A}))$ and $C_{c}(X\ast G,\tilde{r}^{\ast}(\mathcal{K}))$,
where for $f\in C_{cc}(G,r^{\ast}(\mathcal{A})),$ the $L^{I}$-norm of
$f$ is
\[
\max\left\{  \sup_{u\in G^{(0)}}\int\left\|  f(\gamma)\right\|  _{A(u)}%
  \,d\lambda^{u},\sup_{u\in G^{(0)}}\int\left\|
    f^{\ast}(\gamma)\right\|
_{A(u)}\,d\lambda^{u}\right\}  \text{,}%
\]
and similarly for $C_{c}(X\ast G,\tilde{r}^{\ast}(\mathcal{K}))$.
Since every representation of one of these algebras is continuous in
the $L^{I}$-norm by Renault's disintegration theorem
\cite[Theorem~4.1]{ren:jot87}, we conclude that $\Phi$ extends to a
$C^{\ast}$-homomorphism from $C^{\ast}(G;\mathcal{A})$ to
$C^{\ast}(X\ast G;\mathcal{K})$ that is the inverse of $\Psi$.
\end{proof}

\emph{On the basis of Theorem~\ref{Reduction} we may and shall assume
  from now
on that} $X=G^{(0)}$ \emph{and that} $A(u)$ \emph{is an elementary}
  $C^{\ast} 
$\emph{-algebra for every} $u\in G^{(0)}.$

\begin{proof}
  [Proof of the Sufficiency] 
If $G$ acts freely and
properly on $G^{(0)}$, then $G^{(0)}$ is an equivalence between $G$
and the quotient space $G^{(0)}/G$ (viewed as a cotrivial groupoid) in
the sense of \cite{mrw:jot87}. Further, if
$G^{(0)}\ast\mathcal{A}:=\{(u,a)\mid a\in A(u)\}$, then
$G^{(0)}\ast\mathcal{A}$ serves as a bundle of Morita equivalences
between $\mathcal{A}$ and $\mathcal{A}^{G^{(0)}}$ in the sense of
\cite{kmrw:ajm98}. (See \cite{ren:jot87}, too.) Indeed, recall that
$\mathcal{A}% 
^{G^{(0)}}$ is the quotient $(G^{(0)}\ast\mathcal{A})/G$, where
$(r(\gamma
),a)\cdot\gamma=(s(\gamma),\sigma_{\gamma}^{-1}(a))$. 
Then the $\mathcal{A}$-valued inner 
product on $G^{(0)}\ast\mathcal{A}$ is given by the formula
\[
\brip \mathcal{A}<(u,a),(u,b)>=a^{\ast}b\text{,}%
\]
while the $\mathcal{A}^{G^{(0)}}$-valued inner product on $G^{(0)}%
\ast\mathcal{A}$ is given by the formula%
\[
\blip {\mathcal{A}^{G^{(0)}}}<(u,a),(u,b)> 
=[u,ab^{\ast}]\text{,}
\]
where $[u,ab^{\ast}]$ denotes the image of $(u,ab^{\ast})$ in $\mathcal{A}%
^{G^{(0)}}$. By Corollaire 5.4 of \cite{ren:jot87},%
\[
C^{\ast}(G;\mathcal{A})\text{\quad is Morita equivalent to\quad}
C_{0}(G^{(0)}/G,\mathcal{A}%
^{G^{(0)}})\text{.}%
\]
This second algebra is a continuous-trace $C^{\ast}$-algebra since
$\mathcal{A}^{G^{(0)}}$ is a bundle of elementary $C^{\ast}$-algebras
satisfying Fell's condition, as we noted earlier.
\end{proof}

\section{Necessity in Theorem~\ref{Main}}

The proof of the necessity in Theorem~\ref{Main} is modeled on the
proofs in \cite{muhwil:ms90} and \cite{muhwil:ms92}, which in turn are
inspired by ideas in \cite{gre:pjm77}. There are, however, a number of
new difficulties that must be overcome.

For each $u\in G^{(0)}$, we fix an irreducible representation
$\pi_{u}$ of $A(u)$ on a Hilbert space $H_{u}$. Since each $A(u)$ is
elementary the $\pi_{u}$'s are unique up to unitary equivalence. We
define $L^{u}$ to be the representation of $C^{\ast}(G;\mathcal{A})$
on the Hilbert space
$L^{2}(\lambda_{u})\otimes H_{u}$ according to the formula%
\begin{equation}
L^{u}(f)\xi(\gamma)=\int\pi_{u}\circ\sigma_{\gamma}^{-1}(f(\gamma\alpha
))\xi(\alpha^{-1})\,d\lambda^{u}(\alpha)\text{,}\label{LuDef}%
\end{equation}
where $f\in C_{c}(G,r^{\ast}(\mathcal{A}))$ and\ $\xi\in L^{2}(\lambda
_{u})\otimes H_{u}$. (Recall that $\lambda_{u}$ is the image of
$\lambda^{u}$ under inversion.) Thus, $L^{u}$ is the representation of
$C^{\ast }(G;\mathcal{A})$ \emph{induced} by the irreducible
representation $\pi_{u}$ viewed as a representation of
$C_{0}(G^{(0)},\mathcal{A})$. This implies, in particular, that
replacing $\pi_{u}$ by a unitarily equivalent representation does not
affect the unitary equivalence class of $L^{u}$. The following lemma
and corollary capture the salient features of the $L^{u}$ that we
shall use.

\begin{lemma}
\label{L^u}Under the hypothesis that the action of $G$ on $G^{(0)}$ is free
(and Hypothesis~\ref{Hyp1} on $\mathcal{A}$), the following assertions
hold:

\begin{enumerate}
\item Each representation $L^{u}$ is irreducible.
  
\item $L^{u}$ is unitarily equivalent to $L^{v}$ if and only if $u$
  and $v$ lie in the same orbit.
  
\item The map $u\rightarrow L^{u}$ is continuous.
\end{enumerate}
\end{lemma}

\begin{proof}
  The proof follows the lines of the arguments in \cite[Lemma 2.4 and
  Proposition 2.5]{muhwil:ms90}.  Only minor changes need to be made
  to accommodate the presence of $\mathcal{A}$. The key point is that
  $L^{u}$ is unitarily equivalent to the representation $R^{u}$ of
  $C^{\ast}(G;\mathcal{A})$ defined
by the formula%
\[
R^{u}(f)\xi(\gamma\cdot
u)=\int\pi_{u}\circ\sigma_{\gamma}^{-1}(f(\gamma
\alpha))\xi(\alpha^{-1}\cdot u)\,d\lambda^{u}(\alpha)\text{,}%
\]
$f\in C_{c}(G,r^{\ast}(\mathcal{A}))$, $\xi\in L^{2}([u],\mu_{\lbrack
  u]})\otimes H_{u}$, where $[u]$ denotes the orbit of $u$ and
$\mu_{\lbrack u]}$ is the image of $\lambda_{u}$ under the map
$r|s^{-1}(u)$. The fact that the action of $G$ on $G^{(0)}$ is free
(i.e., $G$ is a principal groupoid) implies that $r|s^{-1}(u)$ is a
bijection between $s^{-1}(u)$ and $[u]$. It is a Borel isomorphism, of
course, because of our separability hypotheses and the fact that
$r|s^{-1}(u)$ is continuous.

The value of $R^{u}$ for us lies in the fact that it is evident how to
express $R^{u}$ as the integrated form of a representation of
$(G,\mathcal{A})$ in the sense of \cite[Definition 3.4]{ren:jot87}.
The measure class on $G^{(0)}$ is, of course, that determined by
$\mu_{\lbrack u]}$ and the Hilbert bundle $\mathcal{H}$ is the
constant bundle determined by $H_{u}$ over the orbit of $u$, i.e.,
\[
H(v)=\left\{
\begin{array}
[c]{cc}%
\{v\}\times H_{u} & v\in\lbrack u]\\
0 & \text{otherwise}%
\end{array}
\right.  \text{.}%
\]
Thus, $\int^{\oplus}H(v)\,d\mu_{\lbrack u]}$ is identified with $L^{2}%
([u],\mu_{\lbrack u]})\otimes H_{u}$ in the standard fashion. The
groupoid $G$
is represented on $\mathcal{H}$ according to the formula%
\[
U_{\gamma}((s(\gamma),\xi))=(r(\gamma),\xi)\text{,}%
\]
$\xi\in H_{u}$, $s(\gamma)\in\lbrack u]$, i.e.,
$\{U_{\gamma}\}_{\gamma\in G}$ is just the translation representation,
and $\mathcal{A}$ is represented on
$\mathcal{H}$ according to the formula%
\[
a\cdot(v,\xi)=(v,\pi_{u}\circ\sigma_{\gamma}(a)\xi)\text{,}%
\]
$a\in A(v)$, $(v,\xi)\in H(v)$, where $\gamma$ is the unique element
in $G$ with source $v$ and range $u$.

Observe that the $C_{0}(G^{(0)},\mathcal{A})$ acts as multipliers on
$C^{\ast }(G;\mathcal{A})$ according to the formula $\Phi\cdot
f(\gamma)=\Phi (r(\gamma))f(\gamma)$ for $\Phi\in
C_{0}(G^{(0)},\mathcal{A})$ and $f\in C_{c}(G,r^{\ast}(\mathcal{A}))$.
The extension $\tilde{R}^{u}$ of $R^{u}$ to
the multiplier algebra $C^{\ast}(G;\mathcal{A})$ represents $C_{0}%
(G^{(0)},\mathcal{A})$ on $L^{2}([u],\mu_{\lbrack u]})\otimes H_{u}$
via the
equation%
\[
\tilde{R}^{u}(\Phi)\xi(v)=\pi_{u}\circ\sigma_{\gamma}(\Phi(v))\xi(v)\text{,}%
\]
again, where $\gamma$ is the unique element in $G$ with source $v$ and
range
$u$. It is clear from this that the weak closure of the algebra $\tilde{R}%
^{u}(C_{0}(G^{(0)},\mathcal{A}))$ is the full algebra of decomposable
operators on $L^{2}([u],\mu_{\lbrack u]})\otimes H_{u}$. It follows
that any projection that commutes with
$R^{u}(C^{\ast}(G;\mathcal{A}))$ must be diagonal. On the other hand,
it follows from the definition of the representation of $G$,
$\{U_{\gamma}\}_{\gamma\in G}$, that a diagonal operator commuting
with $R^{u}(C^{\ast}(G;\mathcal{A}))$ must commute with
$\{U_{\gamma}\}_{\gamma\in G}$, and therefore must be constant a.e.\ %
$\mu_{\lbrack u]}$. This proves that $R^{u}$, and hence $L^{u}$, is
irreducible.

If $u$ and $v$ lie in the same orbit, it is clear that translation by
the (unique) $\gamma$ with source $v$ and range $u$ implements an
equivalence between $L^{u}$ and $L^{v}$. On the other hand, if $u$ and
$v$ lie in different orbits, then $L^{u}$ and $L^{v}$ are disjoint.
Indeed, the representations $N_{u}$ and $N_{v}$ of $C_{0}(G^{(0)})$
obtained by restricting $\tilde{R}^{u}$ and $\tilde{R}^{v}$ to
$C_{0}(G^{(0)})$, viewed as a subalgebra of
$M(C^{\ast}(G;\mathcal{A}))$, are supported on the disjoint sets $[u]$
and $[v]$. Arguing just as we did in the proof of \cite[Proposition
2.5]{muhwil:ms90}, using \cite[Lemma 4.15]{wil:tams81}, we conclude
$L^{u}$ and $L^{v}$ are disjoint.

Finally, to see that the map $u\mapsto L^{u}$ is continuous, observe
that Hypothesis~\ref{Hyp1} guarantees that for each point $u\in
G^{(0)}$ we can find a neighborhood $\mathcal{V}_{u}$ of $u$ on which
the $H_{v}$'s can be
chosen to be the fibres of a (topological) Hilbert bundle $\mathcal{\tilde{H}%
}$ and on which we can choose the $\pi_{v}$'s so that for any section
$\Phi\in C_{0}(G^{(0)},\mathcal{A})$ that is supported on
$\mathcal{V}_{u}$ and any two $C_{0}$-sections of
$\mathcal{\tilde{H}}$ over $\mathcal{V}_{u}$, $\xi$ and $\eta$, the
function $v\mapsto \bip(\pi_{v}(\Phi(v))\xi(v)|{\eta(v)})_{H_{v}}$ is
continuous. It follows from the continuity of the Haar system that
given such sections $\xi$ and $\eta$ of $\mathcal{\tilde{H}}$ and any
two functions $g$ and $h$ in $C_{c}(G)$, the function $v\mapsto
\bip(L^{v}(f)(g\otimes
\xi)|{(h\otimes\eta)})$ (where the inner products are taken in $L^{2}%
(\lambda_{v})\otimes H_{v}$) is continuous for all $f\in
C_{c}(G;r^{\ast }(\mathcal{A}))$. This shows that the map $u\mapsto
L^{u}$ is continuous.
\end{proof}

\begin{cor}
\label{L^ubis}Assume that $G$ is principal and that $\mathcal{A}$ is an
elementary $C^{\ast}$-bundle over $G^{(0)}$, satisfying Fell's
condition, on which $G$ acts. If $C^{\ast}(G;\mathcal{A})$ has
continuous trace, then the map that sends $u\in G^{(0)}$ to the
unitary equivalence class of $L^{u}$ defines a continuous open
surjection of $G^{(0)}$ onto $\specnp{C^{\ast }(G;\mathcal{A})}$ that
is constant on $G$-orbits. In particular, orbits are closed and
$G^{(0)}/G$ is homeomorphic to $\specnp{C^{\ast}(G;\mathcal{A})}$.
\end{cor}

\begin{proof}
  The proof is also essentially the same as the proof in
  \cite[Proposition 2.5]{muhwil:ms90}. Here is an outline. Write
  $\Psi$ for the map $u\mapsto \lbrack L^{u}]$. Then by
  Lemma~\ref{L^u}, $\Psi$ is continuous and constant on $G$-orbits.
  Thus $\Psi$ passes to continuous map on $G^{(0)}/G$ with the
  quotient topology (no matter how bad that might be). Since, however,
  $\specnp{C^{\ast}(G;\mathcal{A})}$ is Hausdorff by hypothesis, we
  conclude that $G^{(0)}/G$ is Hausdorff.
  
  Suppose that $L$ is an irreducible representation of
  $C^{\ast}(G;\mathcal{A})$ and let $M$ be the representation of
  $C_{0}(G^{(0)})$ obtained by extending $L$ to the multiplier algebra
  of $C^{\ast}(G;\mathcal{A})$ and then restricting to
  $C_{0}(G^{(0)})$. The kernel $J$ of $M$ is the set of functions in
  $C_{0}(G^{(0)})$ that vanish on a closed set $F$ in $G^{(0)}$. Then
  $F$ is easily seen to be invariant. Indeed, one may do this directly
  or use the fact that it supports the quasi-invariant measure
  associated to the disintegrated form of $L$ guaranteed by
  \cite{ren:jot87}. Further, since $L$ is irreducible, $F$ cannot be
  expressed as the union of two disjoint, closed, $G$-invariant sets.
  Since the quotient map from $G^{(0)}$ to $G^{(0)}/G$ is continuous
  and open, we may apply the lemma on page 222 of \cite{gre:am78} to
  conclude that $F$ is an orbit closure. Since we now know that orbits
  are closed, $F$ is, in fact, an orbit. Thus $L$ factors through
  $C^{\ast}(G|_{[u]};\mathcal{A})$. However, $G|_{[u]}$ is a
  \emph{transitive} principal groupoid. A little reflection, using
  Theorem 3.1 of \cite{mrw:jot87}, reveals that every irreducible
  representation of $C^{\ast}(G|_{[u]};\mathcal{A})$ is unitarily
  equivalent to $L^{u}$.
\end{proof}

We now assume the action of $G$ on $G^{(0)}$ is not proper and use
Lemma 2.6 of \cite{muhwil:ms90} to choose a sequence
$\{\gamma_{n}\}\subseteq G$ such that $\gamma_{n}\rightarrow\infty$ in
the sense that $\{\gamma_{n}\}$ eventually escapes each compact subset
of $G$, and such that $r(\gamma_{n}),s(\gamma _{n})\rightarrow z$, for
some $z\in G^{(0)}$. We shall fix this sequence for the remainder of
the proof.  We also choose a relatively compact neighborhood $U$ of
$z\in G^{(0)}$ and a section $g$ in the Pedersen ideal of
$C_{0}(G^{(0)},\mathcal{A})$ such that $g$ is non-negative, compactly
supported, and satisfies $\tr(\pi_{u}(g(u)))\equiv1$ on $U$. The fact
that $\mathcal{A}$ satisfies Fell's condition guarantees that such
choices are possible.

With these ingredients fixed, we want to build a special neighborhood
$E$ of $z$ in $G$, following the analysis on pages 236--238 of
\cite{muhwil:ms90}. First observe, as we have above, that since $G$ is
principal, $r$ maps $G_{z}$ bijectively onto $[z]$ while $s$ maps
$G^{z}$ bijectively onto $[z]$, where, recall, $[z]$ denotes the orbit
of $z$. Since $[z]$ is closed by Corollary~\ref{L^ubis} while $r$ is
continuous and open on $G$, we see that $r$ maps $G_{z}$
homeomorphically onto $[z]$. Likewise, $s$ maps $G^{z}$
homeomorphically onto $[z]$. Also, since $G$ is principal,
multiplication induces a homeomorphism between $G_{z}\times G^{z}$ and
$G|_{[z]}$. Let $N$ be the closed support of $g$, a compact subset of
$G^{(0)}$, and set $F_{z}:=G_{z}\cap r^{-1}([z]\cap N)$ and
$F^{z}:=G^{z}\cap s^{-1}([z]\cap N)$, obtaining compact subsets of
$G_{z}$ and $G^{z}$, respectively. Then we see that if $\gamma\in
G|_{[z]}$ and if $g(s(\gamma))\neq0$ and $g(r(\gamma ))\neq0$, then
$\gamma\in F_{z}F^{z}$.

According to Lemma 2.7 of \cite{muhwil:ms90}, we may select symmetric,
conditionally compact open neighborhoods $W_{0}$ and $W_{1}$ of
$G^{(0)}$ such that $\overline{W_{0}}\subseteq W_{1}$. (Recall that a
neighborhood $W$ of $G^{(0)}$ is conditionally compact in case $VW$
and $WV$ are relatively compact subsets of $G$ for each relatively
compact subset $V$ in $G$.) We
select such a pair, as we may, with the additional property that $F_{z}%
F^{z}\subseteq W_{0}zW_{0}$. Then from the preceding paragraph, we see
that if $\gamma\notin W_{0}zW_{0}$, then either $g(s(\gamma))=0$ or
$g(r(\gamma))=0$.

By construction,
\[
\overline{W_{1}}^{7}z\backslash W_{0}z\subseteq
r^{-1}(G^{(0)}\backslash
N)\text{.}%
\]
So we may find relatively compact open neighborhoods $V_{0}$ and
$V_{1}$ of $z$ in $G$ so that $V_{0}\subseteq W_{0}$,
$\overline{V_{0}}\subseteq V_{1}$, and so that
\[
\overline{W_{1}}^{7}\overline{V_{1}}\backslash W_{0}V_{0}\subseteq
r^{-1}(G^{(0)}\backslash N)\text{.}%
\]
With these $V_{0}$ and $V_{1}$ so chosen, the special open
neighborhood $E$ of $z$ in $G$ that we want is defined to be
$E:=W_{0}V_{0}W_{0}$.

Observe that we have
\[
\overline{W_{1}}^{7}\overline{V_{1}}\,\overline{W_{1}}^{7}\backslash
E=\overline{W_{1}}^{7}\overline{V_{1}}\,\overline{W_{1}}^{7}\backslash
W_{0}V_{0}W_{0}\subseteq r^{-1}(G^{(0)}\backslash N)\text{.}%
\]
Set%
\[
g^{1}(\gamma):=\left\{
\begin{array}
[c]{cc}%
g(r(\gamma)),\; & \gamma\in\overline{W_{1}}^{7}\overline{V_{1}}\,\overline
{W_{1}}^{7}\\
0,\; & \gamma\notin E
\end{array}
\text{.}\right.
\]
Then, since $\overline{W_{1}}^{7}\overline{V_{1}}\,\overline{W_{1}}%
^{7}\backslash E\subseteq r^{-1}(G^{(0)}\backslash N)$, $g^{1}$ is a
continuous section of $r^{\ast}(\mathcal{A})$ on $G$ that vanishes
outside $E.$

Observe the following containment relations among relatively compact
sets:
$E^{2}=W_{0}V_{0}W_{0}^{2}V_{0}W_{0}\subseteq W_{0}^{4}V_{0}W_{0}^{4}%
\subseteq\overline{W_{0}}^{4}\overline{V_{0}}\,\overline{W_{0}}^{4}\subseteq
W_{1}^{4}V_{1}W_{1}^{4}$. Hence, we may find a compactly supported
function $b$ on $G$ such that $0\leq b(\gamma)\leq1$ for all $\gamma$,
$b\equiv1$ on $E^{2}$ and $b\equiv0$ off $W_{1}^{4}V_{1}W_{1}^{4}$.
Replacing $b$ by $\frac{b+b^{\ast}}{2}$, if necessary, we may assume
that $b$ is a selfadjoint element of the convolution algebra of
scalar-valued functions $C_{c}(G)$.

Define
$F(\gamma):=g(r(\gamma))\sigma_{\gamma}(g(s(\gamma)))b(\gamma)$. By
our choices of $g$ and $b$, $F$ belongs to
$C_{c}(G,r^{\ast}(\mathcal{A}))$, $F$ is selfadjoint and
\[
L^{u}(F)\xi(\gamma)=
\pi_{u}(\sigma_{\gamma}^{-1}(g(r(\gamma))))\int\pi_{u}(\sigma_{\alpha}%
^{-1}(g(r(\alpha))))b(\gamma\alpha^{-1})\xi(\alpha)\,d\lambda_{u}(\alpha)
\]
for all $u$ by the definition of $L^{u}$ (cf. \eqref{LuDef}). Let
$P_{u,1}$ be the projection onto $\mathcal{E}_{u,1}:=$
$L^{2}(G_{u}\cap E)\otimes H_{u}$
and let $P_{u,2}$ be the projection onto the orthocomplement, $\mathcal{E}%
_{u,2}:=$ $L^{2}(G_{u}\backslash E)\otimes H_{u}$. Then, if
$P_{u,1}\xi=\xi
$,$\;$we see that$\;$%
\begin{align*}
  L^{u}(F)\xi(\gamma) &
  =\pi_{u}(\sigma_{\gamma}^{-1}(g(r(\gamma))))\int _{G_{u}\cap
    E}\pi_{u}(\sigma_{\alpha}^{-1}(g(r(\alpha))))\xi(\alpha
  )\,d\lambda_{u}(\alpha)\\
  & =\pi_{u}\circ\sigma_{\gamma}^{-1}(g^{(1)}(\gamma))\int_{G_{u}\cap
    E}\pi
  _{u}\circ\sigma_{\alpha}^{-1}(g^{(1)}(\alpha))\xi(\alpha)d\lambda_{u}(\alpha)
\end{align*}
for all $\gamma\in G_{u}\cap E$ because $b$ is identically $1$ on
$E^{2}$.  However, by definitions of $E$\ and $g^{1}$, the equation
persists when $\gamma\in G_{u}\backslash E$, yielding $0$. Thus,
$P_{u,1}$ commutes with
$L^{u}(F)$. Moreover, when $u=z$, these formulas show that $L^{z}%
(F)P_{z,1}=L^{z}(F)$.

We now want to show that $L^{u}(F)P_{u,1}\geq0$ and we want to analyze
the
trace, $\tr(L^{u}(F)P_{u,1})$. However, when $\xi$ is in the range of $P_{u,1}%
$, the formula for $L^{u}(F)\xi$ shows that
\begin{align*}
  \bip(L^{u}(F)\xi|\xi) & =\int\!\!
\int(\pi_{u}\circ\sigma_{\gamma}^{-1}(g^{(1)}%
(\gamma))\pi_{u}\circ\sigma_{\alpha}^{-1}(g^{(1)}(\alpha))\xi(\alpha
),\xi(\gamma))d\lambda_{u}(\alpha)d\lambda_{u}(\gamma)\\
& =\int\!\!
\int(\pi_{u}\circ\sigma_{\alpha}^{-1}(g^{(1)}(\alpha))\xi(\alpha
),\pi_{u}\circ\sigma_{\gamma}^{-1}(g^{(1)}(\gamma))\xi(\gamma))d\lambda
_{u}(\alpha)d\lambda_{u}(\gamma) \\
&\geq0\text{.}%
\end{align*}
As for the trace, observe that if $K_{u}$ is defined by the formula%
\[
K_{u}(\gamma,\eta)=\pi_{u}(\sigma_{\gamma}^{-1}(g^{1}(\gamma))\sigma_{\eta
}^{-1}(g^{1}(\eta)))
\]
on $G_{u}\times G_{u}$, then our calculations show that $K_{u}$ is
continuous, positive semidefinite and supported on $(G_{u}\cap
E)\times(G_{u}\cap E)$, and that
$\bip(L^{u}(F)P_{u,1}\xi|\zeta)=\int\!\!\!\int(K_{u}(\gamma,\eta)\xi(\eta
),\zeta(\gamma))d\lambda_{u}(\eta)d\lambda_{u}(\gamma)$. Consequently,
we may use Duflo's generalization of Mercer's theorem,
\cite[Proposition 3.1.1]{bcdlrrv:representations}, and the fact that
$K_{u}(\gamma,\gamma)=\pi_{u}\circ
\sigma_{\gamma}^{-1}(g(r(\gamma)))^{2}$ to conclude that$\;\tr(L^{u}%
(F)P_{u,1})=$%
\[
\int_{G_{u}\cap E}\tr(\pi_{u}\circ\sigma_{\gamma}^{-1}(g(r(\gamma
)))^{2})\,d\lambda_{u}(\gamma)\text{.}%
\]
By our choice of $g$, this expression is continuous in $u$ and when
$u=z$ yields the value $\tr(L^{z}(F))$.

We will show that there is a positive number $a$ such that
\begin{equation}
\bigl\| (L^{s(\gamma_{n})}(F)P_{s(\gamma_{n}),2})^{+}\bigr\|
\geq2a \label{norm}%
\end{equation}
eventually, where $(L^{s(\gamma_{n})}(F)P_{s(\gamma_{n}),2})^{+}$
denotes the positive part of the selfadjoint operator
$L^{s(\gamma_{n})}(F)P_{s(\gamma _{n}),2}$ Therefore, eventually
\[
\bigl(\text{the largest eigenvalue of}\enspace L^{s(\gamma_{n}
  )}(F)P_{s(\gamma_{n}),2}\bigr) \geq2a.
\]
Assume that we have shown this, and set
\[
q(t)=\left\{
\begin{array}
[c]{cc}%
0, & t\leq a\\
2(t-a), & a\leq t\leq2a\\
t, & 2a
\end{array}
\right.
\]
Then $q(F)$ is a positive element in the Pedersen ideal of $C^{\ast
}(G,\mathcal{A})$, and so the function $u\rightarrow \tr(L^{u}(q(F)))=$%
\[
\tr(L^{u}(q(F))P_{u,1})+\tr(L^{u}(q(F))P_{u,2})
\]
is finite and continuous in $u$, with value $\tr(L^{z}(q(F)))$ at
$u=z$.
(Recall that $L^{z}(F)=L^{z}(F)P_{z,1}$ and so $L^{z}(q(F))=L^{z}%
(q(F))P_{z,1}$.) On the other hand, we showed that $L^{u}(F)P_{u,1}$
is positive. Since $P_{u,1}$ commutes with $L^{u}(F)$,
$L^{u}(F)P_{u,1}=L^{u}(F^{+})P_{u,1}$. But also we showed that
$u\rightarrow \tr(L^{u}(F)P_{u,1})$ is continuous at $z$.
Consequently, so is $u\rightarrow \tr(L^{u}(F^{+})P_{u,1})$. Since
$q(F)\leq F^{+}$, the function $u\rightarrow \tr(L^{u}(q(F))P_{u,1})$
is continuous by Lemma 4.4.2(i) in \cite{dix:cs-algebras}, with value
$\tr(L^{z}(q(F)))=\tr(L^{z}(q(F))P_{z,1})$ at $u=z$, also. Therefore
\[
\lim_{u\rightarrow z}\tr(L^{u}(q(F))P_{u,2})=0\text{.}%
\]
Since the largest eigenvalue of $L^{s(\gamma_{n})}(F)P_{s(\gamma_{n}),2}%
\geq2a$, the largest eigenvalue of $L^{s(\gamma_{n})}(q(F))P_{s(\gamma_{n}%
  ),2}\geq2a$ also. Consequently,
\[
\liminf_{n}\tr(L^{s(\gamma_{n})}(q(F))P_{s(\gamma_{n}),2})\geq2a\text{.}%
\]
This contradiction will complete the proof.

So we will finish by verifying the asserted inequality (\ref{norm}).
To this end, choose an open neighborhood of $z$ in $G$, $V_{2}$, that
is contained in $V_{0}$ and choose a conditionally compact
neighborhood $Y$ of $G^{(0)}$ such that if $v\in V_{2}$, then $r$ maps
$Yv$ into $U$. Without loss of generality,
we may assume that $Y\subseteq W_{0}$. Observe that if $\gamma_{n}%
\notin\overline{W_{1}}^{2}\overline{V_{1}}\,\overline{W_{1}}^{2}$ then
for $\gamma\in Y\gamma_{n}$, $\gamma\notin E$. Indeed, if
$\gamma=\gamma^{\prime
}\gamma_{n}\in E\cap Y\gamma_{n}$, then $\gamma_{n}\in(\gamma^{\prime}%
)^{-1}E\subseteq W_{0}^{2}V_{0}W_{0}\subseteq\overline{W_{1}%
}^{2}\overline{V_{1}}\,\overline{W_{1}}^{2}$ contrary to assumption.
So, since $r(\gamma_{n})$ and $s(\gamma_{n})$ are tending to $z$,
while $\gamma_{n}$ eventually escapes
$\overline{W_{1}}^{2}\overline{V_{1}}\,\overline{W_{1}}^{2}$, we can
conclude that for $n$ sufficiently large, whenever $\gamma$ lies in
$Y\gamma_{n}$, then $\gamma\notin E$ while $r(\gamma)$ and $s(\gamma)$
lie in $U$. From now on, we will assume that $n$ is sufficiently large
so that these conditions are satisfied.

Next observe that since for each $n$, the map $\gamma\rightarrow\pi
_{s(\gamma_{n})}\circ\sigma_{\gamma}^{-1}(g(r(\gamma))$ defines a
continuous family of rank $1$ projections on the Hilbert space
$H_{s(\gamma_{n})}$, we can find a Borel family of unit vectors
$\gamma\rightarrow v_{\gamma}^{n}$ such that
$\pi_{s(\gamma_{n})}\circ\sigma_{\gamma}^{-1}(g(r(\gamma))$ is the
rank $1$ projection determined by $v_{\gamma}^{n}$.

Let $h_{n}(\gamma)=1_{Y\gamma_{n}}(\gamma)\times v_{\gamma}^{n}$ where
$1_{Y\gamma_{n}}$ denotes the characteristic function of
$Y\gamma_{n}$.  Observe that if $\gamma,\alpha\in Y\gamma_{n}$, then
$\gamma\alpha^{-1}\in Y\gamma_{n}\gamma_{n}^{-1}Y\subseteq
YV_{0}Y\subseteq W_{0}V_{0}W_{0}=E$.  Consequently,
$b(\gamma\alpha^{-1})=1$ and we may calculate to find that if
$\gamma\in Y\gamma_{n}$, then%
\begin{align*}
  L^{s(\gamma_{n})}&(F)h_{n} (\gamma)\\
&=\pi_{s(\gamma_{n})}(\sigma_{\gamma}%
^{-1}(g(r(\gamma))))
%&\qquad\qquad
 \int\pi_{s(\gamma_{n})}(\sigma_{\alpha}^{-1}%
(g(r(\alpha))))b(\gamma\alpha^{-1})h_{n}(\alpha)\,d\lambda_{s(\gamma_{n}%
  )}(\alpha)\\
&  =v_{\gamma}^{n}\lambda_{r(\gamma_{n})}(Y)\text{.}%
\end{align*}
Hence $\bip( L^{s(\gamma_{n})}(F)h_{n}|h_{n}) =\lambda_{r(\gamma
  _{n})}(Y)^{2}$.  However, by our assumption on $\gamma_{n}$, $h_{n}$
lies in $\mathcal{E}_{s(\gamma_{n}),2}$, so $\bip(
(L^{s(\gamma_{n})}(F)P_{s(\gamma _{n}),2})^{+}h_{n}|h_{n}) \geq\bip(
L^{s(\gamma_{n})}(F)P_{s(\gamma _{n}),2}h_{n}|h_{n}) =\bip(
L^{s(\gamma_{n})}(F)h_{n}|h_{n}) =\lambda_{r(\gamma_{n})}(Y)^{2}$.
This shows that $\bigl\|
(L^{s(\gamma_{n})}(F)P_{s(\gamma_{n}),2})^{+}\bigr\| \geq
\lambda_{r(\gamma_{n})}(Y)$ provided $n$ is sufficiently large. But
the continuity of the Haar system implies that
$\liminf\lambda_{r(\gamma_{n} )}(Y)>0$, as $n\rightarrow\infty$. This
verifies equation (\ref{norm}) and completes the proof of
Theorem~\ref{Main}.
%TCIMACRO{\TeXButton{End Proof}{\endproof}}%
%BeginExpansion
%\endproof
%EndExpansion 

%\bibliographystyle{amsalpha}      % For drafts
%\bibliographystyle{amsplain}      % For final copy
%\bibliography{references-nov01}

\def\noopsort#1{}\def\cprime{$'$}
\providecommand{\bysame}{\leavevmode\hbox to3em{\hrulefill}\thinspace}
\providecommand{\MR}{\relax\ifhmode\unskip\space\fi MR }
\providecommand{\MRhref}[2]{%
  \href{http://www.ams.org/mathscinet-getitem?mr=#1}{#2}
}
\providecommand{\href}[2]{#2}

\end{document}